\documentclass[12pt]{article}
\usepackage{amsmath}
\usepackage{amssymb}
\usepackage{amsgen}
\usepackage{amsfonts}

\newtheorem{theorem}{Theorem} 
\newtheorem{corollary}[theorem]{Corollary} 
\newtheorem{conjecture}[theorem]{Conjecture} 
\newtheorem{lemma}[theorem]{Lemma} 
\newtheorem{proposition}[theorem]{Proposition} 
 
\newtheorem{example}[theorem]{Example}

\newtheorem{remark}{Remark} 

\numberwithin{equation}{section}

\def\ppp{{\mathbb{P}}}

\def\fff{{\mathbb{F}}}

\def\ccc{\mathbb{C}}
\def\zzz{\mathbb{Z}}

\def\pf{\noindent {\bf Proof}:\ }
\def\qed{$\Box$}

\begin{document}

\author{David Joyner, Will Traves\thanks{The first 
author was partially supported by an NSA-MSP grant. 
The second author was supported by a USNA-NARC grant.
Mathematics Dept, US Naval Academy, Annapolis, MD 21402,
wdj@usna.edu, traves@usna.edu. This paper is a significant revision of
the first version, dated early 2003.}}

\title{Representations of finite groups on Riemann-Roch
spaces}

\maketitle

\begin{abstract}
We study the action of a finite group on the Riemann-Roch
space of certain divisors on a curve. 
If $G$ is a finite subgroup of the
automorphism group of a projective curve $X$ 
over an algebraically closed field and $D$ is a 
divisor on $X$ left stable by $G$ then 
we show the irreducible constituents of the natural representation
of $G$ on the Riemann-Roch space $L(D)=L_X(D)$ are of 
dimension $\leq d$, where
$d$ is the size of the smallest $G$-orbit acting
on $X$. We give an example to show that this is, in
general, sharp
(i.e., that dimension $d$ irreducible constituents can occur).
Connections with coding theory, in particular 
to permutation decoding of AG codes, are 
discussed in the last section.
Many examples are included.
\end{abstract}

\tableofcontents

\vskip .2in

Let $X$ be a smooth projective (irreducible) curve over an algebraically
closed field $F$
and let $G$ be a finite subgroup of automorphisms of $X$
over $F$. We often identify $X$ with its set of
$F$-rational points $X(F)$.
If $D$ is a divisor of $X$ which $G$ leaves 
stable then $G$ acts on the Riemann-Roch space $L(D)$.
We ask the question: which (modular) representations 
arise in this way? 

Similar questions have been investigated previously.
For example, 
the action of $G$ on the space of regular differentials,
$\Omega^1(X)$ (which is isomorphic to $L(K)$,
where $K$ is a canonical divisor).
This was first looked at from the representation-theoretic
point-of-view by Hurwitz (in the case $G$ is cyclic) and
Weil-Chevalley (in general).
They were studying monodromy representations on
compact Riemann surfaces. For more details and further
references, see the book by Breuer \cite{B} and the paper
\cite{MP}. Other related works,
include those by Nakajima \cite{N}, Kani \cite{Ka}, K\"ock
\cite{K}, and Borne \cite{Bo1}, \cite{Bo2}, \cite{Bo3}.

The motivation for our study lies in coding theory. 
The construction of AG codes uses the 
Riemann-Roch space $L(D)$ associated to a divisor $D$ of
a curve $X$ defined over a finite field \cite{G}.
Typically $X$ has no non-trivial 
automorphisms\footnote{Indeed, a theorem of Rauch, Popp, and Oort
(see \S 1.2 in \cite{Bo3}, for example)
implies that if $g>3$ 
then the singular points of the moduli space $M_g$
of curves of genus $g$ correspond to curves 
having a non-trivial
automorphism group.}, but when it does we may ask
how this can be used to better understand AG codes 
constructed from $X$. If $G$ is a finite group
acting transitively on a basis of $L(D)$ 
(admittedly an optimistic expectation, but one which gets
the idea across) then one might expect that
fast encoding and decoding algorithms exists for the associated
AG codes. Of course, for such an application, one wants $F$ to be 
finite (and not algebraically closed). 
These ideas are discussed in \S \ref{sec:4} below for
AG codes constructed from the hyperelliptic curves
$y^2=x^p-x$ over $GF(p)$.
Several conjectures on the complexity of permutation decoding of the 
associated AG codes are given there.

\section{The action of $G$ on $L(D)$}
\label{sec:1}

Let $X$ be a smooth projective curve over an algebraically closed
field $F$. Let $F(X)$ denote the function field of
$X$ (the field of rational functions on $X$) and, if
$D$ is any divisor on $X$ then the Riemann-Roch space
$L(D)$ is a finite dimensional $F$-vector space given by
\[
L(D)=L_X(D)= \{f\in F(X)^\times \ |\ div(f)+D\geq 0\}\cup \{0\},
\]
where
$div(f)$ denotes the (principal) divisor of the function
$f\in F(X)$. Let $\ell(D)$ denote its dimension.
We recall the Riemann-Roch theorem,
\[
\ell(D)-\ell(K-D)={\rm deg}(D)+1-g,
\]
where $K$ denotes a canonical divisor and $g$ the 
genus\footnote{We often also use $g$ to denote an element of an
automorphism group $G$. Hopefully, the context will
make our meaning clear.}.

The action of ${\rm Aut}(X)$ on $F(X)$ is 
defined by 
\[
\begin{array}{cccc}
\rho:&{\rm Aut}(X)&\longrightarrow &{\rm Aut}(F(X)),\\
 & g &\longmapsto & (f\longmapsto f^g)
\end{array}
\]
where $f^g(x)=(\rho(g)(f))(x)=f(g^{-1}(x))$.

Note that $Y=X/G$ is also smooth and $F(X)^G=F(Y)$.

Of course, ${\rm Aut}(X)$ also acts on the group $Div(X)$ of
divisors of $X$, denoted $g(\sum_P d_P P)=\sum_P d_Pg(P)$,
for $g\in {\rm Aut}(X)$, $P$ a prime divisor, and $d_P \in \zzz$.
It is easy to show that $div(f^g)=g(div(f))$.
Because of this, if $div(f)+D\geq 0$
then $div(f^g)+g(D)\geq 0$, for all $g\in {\rm Aut}(X)$.
In particular, if the action of $G\subset {\rm Aut}(X)$ on $X$ leaves 
$D\in Div(X)$ stable then $G$ also acts on $L(D)$.
We denote this action by
\[
\rho:G\rightarrow {\rm Aut}(L(D)).
\]

\section{Examples and special cases}

Before tackling the general case, we study the Riemann-Roch 
representations of $G$ when $X=\ppp^1$ or $D$ is the 
canonical divisor.

\subsection{The canonical embedding}

This case was solved by Weil and Chevalley - see the beautiful 
discussion in \cite{MP}.

Let $K$ denote a canonical
divisor of $X$, so deg$(K)=2g-2$ and dim$(L(K))=g$.
Let $\{\kappa_1,...,\kappa_g\}$ denote a basis for $L(K)$.
If the genus $g$ of $X$ is at least $2$ then
the morphism
\[
\begin{array}{ccc}
\phi:X&\longrightarrow &\ppp(\Omega^1(X))\cong \ppp^{g-1}\\
x & \longmapsto & (\kappa_1(x):...:\kappa_g(x))
\end{array}
\]
defines an embedding, the ``canonical embedding'', and
$\phi$ is called the ``canonical map''.
It is known that $L(K)$ is isomorphic (as $F$-vector spaces)
to the space $\Omega^1(X)$ of regular Weil differentials on $X$.
This is contained in the space of all Weil differentials,
$\Omega(X)$. (In the notation of \cite{Sti},
$\Omega^1(X)=\Omega(X)(0)$.)
Since $G$ acts on the set of places of $F$, it acts on the
adele ring of $F$, hence on the space $\Omega(X)$.

Now, even though $K$ might not be fixed by $G$,
there is an action of $G$ on $L(K)$
obtained by pulling back the action of $G$ on $\Omega^1(X)$
via an isomorphism $L(K)\cong \Omega^1(X)$.

The group ${\rm Aut}(X)$ acts on $X$ and on its
image $Y=\phi(X)$ under an embedding $\phi:X\rightarrow \ppp^n $. 
If $\phi$ arises from a very ample linear system 
then an 
automorphism of $Y$ may be represented (via the linear system) 
by an element of $PGL(n+1,F)$ acting
on $\ppp^n$ which preserves $Y$.
For instance, if $D$ is any divisor with deg$(D)>2g$ then
the morphism

\[
\begin{array}{ccc}
\phi:X&\rightarrow &\ppp^{n-1}\\
x & \longmapsto & (f_1(x):...:f_n(x))
\end{array}
\]
defines an embedding,
where $\{f_1,...,f_n\}$ is a basis for $L(D)$
(see, for example, Stepanov \cite{St}, \S 4.4).
This projective representation of $G$ on $L(D)$ 
exists independent of whether or not $D$ is left
stable by $G$. 

\begin{example}
Let $X=\ppp^1/\ccc$ have projective coordinates
$[x:y]$, let $G=\{1,g\}$, where $g(x/y)=y/x$, and
let $D=2[1:0]-[0:1]$, so
$L(D)$ has basis $\{x/y,x^2/y^2\}$.
Then $g(x/y)=(y/x)^3(x^2/y^2)$ and 
$g(x^2/y^2)=(y/x)^3(x/y)$. 
Thus, as an element of
$PGL(2,\ccc)$, $g$ is 
$\left(
\begin{array}{cc}
0 & 1\\
1 & 0
\end{array}
\right)$.

\end{example}

Suppose, for example, $X$ is non-hyperelliptic of genus $\geq 3$
and $\phi$ arises from the canonical embedding.
In this case, we have (a) the
projective representation 
\[
\pi:G\rightarrow {\rm Aut}(\ppp(\Omega^1(X)))
\]
(acting on the canonical embeding of $X$) and 
(b) the projective representation 
obtained by composing the ``natural'' representation
$G\rightarrow {\rm Aut}(\Omega^1(X))$ with the quotient map
${\rm Aut}(\Omega^1(X))\rightarrow {\rm Aut}(\Omega^1(X)/F^\times)
={\rm Aut}(\ppp(\Omega^1(X)))$. These two representations 
are the same.

\begin{remark}
For further details on the representation
$G\rightarrow {\rm Aut}(\Omega^1(X))$, see for example,
the Corollary to Theorem 2 in \cite{K}, Theorem 2.3 in
\cite{MP}, and the book by T. Breuer \cite{B}.

\end{remark}

\subsection{The projective line}

Before tackling the general case, we study the Riemann-Roch 
representations of $G$ when $X=\ppp^1$.

Let $X=\ppp^1/F$, so ${\rm Aut}(X)=PGL(2,F)$, where 
$F$ is algebraically closed. 
Let $\infty=[1:0]\in X$ denote the element corresponding to 
the localization $F[x]_{(1/x)}$.
In this case, the canonical divisor 
is given by $K=-2\infty$, so the Riemann-Roch theorem 
becomes
\[
\ell(D)-\ell(-2\infty -D)={\rm deg}(D)+1.
\]
It is known (and easy to show) that if ${\rm deg}(D)<0$ then 
$\ell(D)=0$ and if deg$(D)\geq 0$ then $\ell(D)={\rm deg}(D)+1$.

In the case of the projective line, there is another 
way to see the action $\rho$ of ${\rm Aut}(X)$ on $F(X)$. 
Each function
$f\in F(X)$ may be written uniquely as a rational
function $f(x)=p(x)/q(x)$, where $p(x)$ and $q(x)$
are polynomials that factor as the product of
linear polynomials. Assume that both $p$ and $q$ are monic, and assume that
the linear factors of them are as well.
The group ${\rm Aut}(X)$ ``acts'' on the set of such functions
$f$ by permuting its zeros and poles according to the action of $G$ on
$X$. (We leave aside how $G$ acts on the constants, so
this ``action'' is not linear.) 
We call this the ``permutation action'',
$\pi: g \longmapsto \pi(g)(f)= f_g$,
where $f_g(x)$ denotes the function $f$ with
zeros and poles permuted by $g$. 

\begin{lemma}
\label{lemma:permaction}
If $G\subset {\rm Aut}(X)$ leaves $D\in Div(X)$ stable 
then
\[
\pi(g)(f)=c\rho(g)(f),
\]
for some constant $c$.
\end{lemma}

\pf
Note that, by definition,
$div(\pi(g)(f))=div(f_g)=g(div(f))$, for $g\in G$ and $f\in L(D)$.
Since $div(\pi(g)(f))=g(div(f))=div(\rho(g)(f))=div(f^g)$, the functions
$f^g$ and $f_g$ must differ by a constant factor.
\qed

\medskip

The above lemma is useful since it is easier to deal with 
$\pi$ than $\rho$ in this case.

A basis for the Riemann-Roch space is explicitly known
for $\ppp^1$. For notational simplicity, let
\[
m_P(x)=
\left\{
\begin{array}{cc}
x, & P=[1:0]=\infty,\\
(x-p)^{-1}, & P=[p:1].
\end{array}
\right.
\]

\begin{lemma}
\label{lemma:P1basis}
Let $P_0=\infty=[1:0]\in X$ denote the point corresponding to 
the localization $F[x]_{(1/x)}$.
For $1\leq i\leq s$, let
$P_i=[p_i:1]$ denote the point corresponding to 
the localization $F[x]_{(x-p_i)}$, for $p_i\in F$.
Let $D=\sum_{i=0}^s a_i P_i$ be a 
divisor, $a_k\in\zzz$ for $0\leq k\leq s$.
\begin{itemize}
\item[(a)]
If $D$ is effective then
\[
\{
1,m_{P_i}(x)^k\ |\ 1\leq k\leq a_i,
0\leq i\leq s\}
\]
is a basis for $L(D)$.

\item[(b)]
If $D$ is not effective but deg$(D)\geq 0$ then write
$D=dP+D'$, where deg$(D')=0$, $d>0$, and $P$ is any point. Let 
$q(x)\in L(D')$ (which is a 1-dimensional vector space) be any 
non-zero element.
Then
\[
\{
m_P(x)^{i}q(x)\ |\ 0\leq i\leq d\}
\]
is a basis for $L(D)$.

\item[(c)]
If deg$(D)<0$ then $L(D)=\{0\}$.
\end{itemize}

\end{lemma}

The first part is Lemma 2.4 in \cite{L}. The other parts follow from the
definitions and the Riemann-Roch theorem.

In general, we have the following result.

\begin{theorem}
\label{thrm:P1}
Let $X$, $F$, $G\subset {\rm Aut}(X)=PGL(2,F)$, and 
$D=\sum_{i=0}^s a_i P_i$ be a divisor as above. 
Let $\rho:G\rightarrow {\rm Aut}(L(D))$ denote the
associated representation. This acts trivially on
the constants (if any) in $L(D)$; we denote this action by $1$.
Let $S={\rm supp}(D)$ and let
\[
S=S_1\cup S_2\cup ... \cup S_m
\]
be the decomposition of $S$ into primitive $G$-sets.

\begin{itemize}
\item[(a)]
If $D$ is effective then 
\[
\rho\cong 1\oplus_{i=1}^m \rho_i,
\]
where $\rho_i$ is a monomial representation
on the subspace 
\[
V_i=\langle m_P(x)^{\ell_j}\ |\ 1\leq \ell_j\leq a_j,
\ P\in S_i\rangle,
\]
satisfying dim$(V_i)=\sum_{P_j\in S_i} a_j$, for $1\leq i \leq m$.
Here $\langle ... \rangle$ denotes the vector space span.

\item[(b)]
If deg$(D)>0$ but $D$ is not effective 
then $\rho$ is a subrepresentation of 
$\rho:G\rightarrow Aut_F\, L(D')$, where $D'$ is a $G$-equivariant 
effective divisor satisfying $D'\geq D$.
\end{itemize}

\end{theorem}

\pf
(a)
Fix an $i$ such that $1\leq i \leq m$.
Consider the subspace $V_i$ of $L(D)$.
Since $G$ acts by permuting the points in $S_i$
transitively, this action induces an action $\rho_i$ on 
$V_i$. This action on $V_i$ is a monomial representation,
by Lemma \ref{lemma:permaction}. It is irreducible since
the action on $S_i$ is transitive, by definition.
Clearly $\oplus_{i=1}^m \rho_m$ is a subrepresentation
of $\rho$. For dimension reasons, this subrepresentation
must be all of
$\rho$, modulo the constants (the trivial representation).

(b)
Since $D$ is not effective,
we may write $D=D^+-D^-$, where $D^+$ and $D^-$ are 
non-zero effective divisors. 
The action of $G$ must preserve $D^+$ and $D^-$. 
Since $L(D)$ is a $G$-submodule of $L(D^+)$,
the claim follows.
\qed

\section{The general case}

Let $X$ be a smooth projective curve defined over
a field $F$. 
The following is our most general result.

\begin{theorem}
\label{thrm:main}
Suppose $G\subset {\rm Aut}(X)$
is a finite subgroup, and that the 
divisor $D\not= 0$ on $X$ is stable under $G$. 
Let $d_0$ denote the size of a smallest $G$-orbit in $X$. 
Each irreducible composition factor of the
representation of $G$ on
$L(D)$ has dimension $\leq d_0$.
\end{theorem}

\begin{remark}
This is best possible in the sense that irreducible subspaces 
of dimension $d_0$ can occur, by Theorem \ref{thrm:P1}
(see also \S \ref{example:shimura} and Example \ref{sec:gfp2}
below).
\end{remark}

\begin{remark}
If $F$ has characteristic $0$ then
every finite dimensional representation of a finite
group is semi-simple (Prop 9, ch 6, \cite{Se}).
If $F$ has characteristic $p$ and $p$ does not divide $|G|$
then every finite dimensional representation of $G$
is semi-simple (Maschke's Theorem, Thrm 3.14, \cite{CR}, or
\cite{Se}, \S 15.7).

\end{remark}

\pf
Let $D_0\not= 0$ be a
effective $G$-invariant divisor of minimal
degree $d_0$.
Let $d=[{\rm deg}(D)/d_0]$ denote the integer part.
The group $G$ acts on each space in the 
composition series

\[
\begin{array}{c}
\{0\}=L(-(d+1)D_0+D)
\subset L(-dD_0+D)
\subset 
L(-(d-1)D_0+D)\subset \\
...\subset 
L(-(d-m)D_0+D)
\subset ...
\subset L(D)\ .
\end{array}
\]
In particular, $G$ acts on the successive quotient spaces

\[
L(-(d-m-1)D_0+D)/L(-(d-m)D_0+D),\ \ \ \ 0\leq m\leq d-1,
\]
by the quotient representation. 
These are all of dimension at most $d_0$
(Prop. 3, ch 8, \cite{F}).

\qed

\begin{corollary}
\label{cor:inequality}
Suppose that $G$ is a non-abelian group acting on a smooth 
projective curve $X$ defined over an algebraically closed
field $F$ and assume $p$ does not divide the order of 
$G$. Let $d_0$ be as in Theorem \ref{thrm:main} and let
$d_G$ denote the largest degree of all irreducible 
($F$-modular) representations of $G$. Then
\[
d_0\geq d_G.
\]
\end{corollary}

\pf
Construct an effective divisor $D$ of $X$ fixed by $G$.
We may assume that its degree is so large that
the formula of Borne \cite{Bo1} implies that each
irreducible representations of $G$ occurs at least
once in the decomposition of $L(D)$. Therefore the set of irreducible
subrepresentations of $L(D)$ are the same as the
set of irreducible representations of $G$.
The result now follows from our theorem.
\qed

\begin{remark}
There are more general conditions for which
$d_0= d_G$ holds. 
For example, assume that $P$ is a point in an orbit of size $d_0$
and let $H=G_P$ denote the stabilizer of $P$, so $d_0=|G|/|H|$. 
Let $\sigma$ denote an irreducible representation of $H$.
If $H$ is normal in $G$ with cyclic quotient and if
all the equivalence classes $\sigma^g$ ($g\in G/H$) are distinct then
$Ind_H^G\, \sigma$ is irreducible and of dimension $d_0$,
by Clifford's theorem. If there is an irreducible
representation of $G$ of this form $Ind_H^G\, \sigma$ 
under the above conditions then $d_0\leq d_G$.
\end{remark}

{\bf Question}: Is there an analog of Corollary \ref{cor:inequality}
for wildly ramified $\pi:X\rightarrow X/G$?

\subsection{Examples}

\begin{example}
\label{example:fermat}
Let $\fff$ be a separable algebraic closure of $\fff_3$. 
Let $X$ denote the Fermat curve over $\fff$ 
whose projective model is given by
$x^4+y^4+z^4=0$. The point
$P=(1:1:1)\in X(\fff)$ is fixed by the action of $G=S_3$.
 
Based on the Brauer character table of $S_3$ over $\fff_3$
(available in GAP \cite{GAP4}), the group $G$ 
has no $2$-dimensional irreducible (modular) representations.
Consequently, $d_G=d_0=1$.
\end{example}

\begin{example}
\label{example:shimura}
Let $k=\ccc$ denote the complex field and let
$X(N)$ denote the modular curve associated to the
principal congruence group $\Gamma(N)$ (see for example
Stepanov, \cite{St}, chapter 8).
It is well-known that the 
group $PSL(2,\zzz/N\zzz)$ is contained in the
automorphism group of $X(N)$. 
Let $X=X(p)$, where $p\geq 7$ is a prime,
and let $G=PSL(2,\fff_p)$. 
In this case, we have, in the notation of the
above corollary, $d_G=p+1$.
(The representations of this simple group
are described, for example, in 
Fulton and Harris \cite{FH} 
\footnote{Actually those of $SL(2,\fff_p)$
are described in \cite{FH}, but it is easy to determine the
representations of $PSL(2,\fff_p)$ from those of $SL(2,\fff_p)$.}.)
\end{example}

\subsection{$y^2=x^p-x$}

In general, if $X$ is a curve defined over a field $F$
with finite automorphism group $G=Aut_F(X)$ then we call
$G$ {\bf large} if $|G|>|X(F)|$.

\begin{lemma}
\label{lemma:large}
If $G$ is large then every point of $X(F)$
is ramified for the covering $X\rightarrow X/G$.
\end{lemma}

\pf
Suppose $P\in X(F)$ is not ramified, so the stabilizer of $P$,
$G_P$, is trivial. In this case,
$|G\cdot P|= |G|/|G_P|= |G|$. But $G\cdot P\subset X(F)$
so $|G\cdot P|\leq |X(F)|$, a contradiction. \qed

\subsubsection{Case $F=GF(p)$}
\label{sec:F=GF(p)}

Let $p\geq 3$ be a prime,
$F=GF(p)$, and let $X$ denote the curve defined by

\[
y^2=x^p-x.
\]
This has genus $\frac{p-1}{2}$. 
We assume that the automorphism group $G=Aut_F(X)$ is
a central 2-fold cover of $PSL(2,p)$,
we have a short exact 
sequence,

\begin{equation}
\label{eqn:*}
1\rightarrow Z \rightarrow G\rightarrow PSL_2(p)\rightarrow 1,
\end{equation}
where $Z$ denotes the center of $G$ 
($Z$ is generated by the hyperelliptic involution).
The following transformations are elements of
$G$:

\begin{equation}
\label{eqn:3.2}
\begin{array}{cc}
\gamma_1=
\left\{
\begin{array}{c}
x\longmapsto x,\\
y\longmapsto -y,
\end{array}
\right. ,
&
\gamma_2=\gamma_2(a)=
\left\{
\begin{array}{c}
x\longmapsto a^2x,\\
y\longmapsto ay,
\end{array}
\right.\\
\gamma_3=
\left\{
\begin{array}{c}
x\longmapsto x+1,\\
y\longmapsto y,
\end{array}
\right. ,
&
\gamma_4=
\left\{
\begin{array}{c}
x\longmapsto -1/x,\\
y\longmapsto y/x^{\frac{p+1}{2}},
\end{array}
\right. 
\end{array}
\end{equation}
where $a\in F^\times$ is a primitive $(p-1)-st$ root of unity.
This group acts transitively on $X(F)$,
so it has an orbit of size $d_0=|X(F)|=p+1$.

Let $P_1=(1:0:1)$ and let $H$ be its stabilizer in $G$.
A counting argument shows that $H$ is a solvable group of order 
$2p(p-1)$ generated by 
$\gamma_1$, $\gamma_2(a)$ and $\gamma_3$.
By Lemma \ref{lemma:large}, every point in 

\[
X(F)=\{(1:0:0),(0:0:1),(1:0:1),...,(p-1:0:1)  \}
\]
is ramified over the covering $X\rightarrow X/G$
in the sense that each stabilizer $G_P=Stab_G(P)$
is non-trivial, $P\in X(F)$. 

It is known (Proposition VI.4.1, \cite{Sti}) that, 
for each $m\geq 1$, the Riemann-Roch space
of $D=mP_1$ has a basis consisting of monomials,

\[
x^iy^j,\ \ \ \ \ \ 
0\leq i\leq p-1,\ j\geq 0,\ 2i+pj\leq m.
\]

\begin{lemma}
The semisimplification $\rho_{ss}$ of 
the representation $\rho$ of $H$ acting on $L(D)$ 
is the direct sum of 
one-dimensional representations of $G$.
\end{lemma}

\pf The generator $\gamma_1$ acts trivially on the 
basis of $L(D)$, whereas
\[
\gamma_2(a):
\left(
\begin{array}{c}
1\\
x\\
\vdots \\
x^ry^s
\end{array}
\right)
\longmapsto 
\left(
\begin{array}{c}
1\\
a^2x\\
\vdots \\
a^{2r+s}x^ry^s
\end{array}
\right)
=
\left(
\begin{array}{cccc}
1 & 0 & ...& 0\\
0 & a^2 & ...& \vdots\\
\vdots &\ddots  & &0 \\
0 & ... & 0 & a^{2r+s}
\end{array}
\right)
\left(
\begin{array}{c}
1\\
x\\
\vdots \\
x^ry^s
\end{array}
\right),
\]
and

\[
\gamma_3:
\left(
\begin{array}{c}
1\\
x\\
\vdots \\
x^ry^s
\end{array}
\right)
\longmapsto 
\left(
\begin{array}{c}
1\\
x+1\\
\vdots \\
(x+1)^ry^s
\end{array}
\right)
=
\left(
\begin{array}{cccc}
1 & 0 & ...& 0\\
1 & 1 & ...& \vdots\\
\vdots &\ddots  & &0 \\
0 & ... & r & 1
\end{array}
\right)
\left(
\begin{array}{c}
1\\
x\\
\vdots \\
x^ry^s
\end{array}
\right),
\]
where the non-zero terms in bottom row of the matrix representation of
$\gamma_3$ are in the last $r+1$ row entries and
consist of the binomial coefficients $\frac{r!}{(r-j)!j!}$,
$0\leq j\leq r$.
Therefore, the group generated by these matrices is 
lower-triangular, hence solvable. \qed

\subsubsection{Case $F=GF(p^2)$}
\label{sec:gfp2}

Let $F=GF(p^2)$ and let $F_0=GF(p)$.

The automorphism group $G=Aut_F(X)$ is
a central 2-fold cover of $PGL(2,p)$ andwe have a short exact sequence,

\begin{equation}
\label{eqn:**}
1\rightarrow Z \rightarrow G
\stackrel{\tau}{\rightarrow} PGL_2(p)\rightarrow 1,
\end{equation}
where $Z$ denotes the subgroup of $G$ generated by the hyperelliptic
involution (which coincides with the center of $G$), 
by G\"ob \cite{G}.
The group $G$ has order $2|PGL(2,p)|=2p(p^2-1)$.
The following transformations generate $G$:

\[
\begin{array}{cc}
\gamma_1=
\left\{
\begin{array}{c}
x\longmapsto x,\\
y\longmapsto -y,
\end{array}
\right. ,
&
\gamma_2=\gamma_2(a)=
\left\{
\begin{array}{c}
x\longmapsto a^2x,\\
y\longmapsto ay,
\end{array}
\right.\\
\gamma_3=
\left\{
\begin{array}{c}
x\longmapsto x+1,\\
y\longmapsto y,
\end{array}
\right. ,
&
\gamma_4=
\left\{
\begin{array}{c}
x\longmapsto -1/x,\\
y\longmapsto y/x^{\frac{p+1}{2}},
\end{array}
\right. 
\end{array}
\]
where $a\in F^\times$ is a primitive $2(p-1)-st$ root
of unity.

\begin{proposition}
\label{prop:orbits} 
Let $p>3$ be a prime.
\begin{itemize}
\item[(a)] Case $p\equiv 3\, ({\rm mod}\, 4)$:

Let $P_1=(1:0:1)$ and fix some $P_2\in X(F)-X(F_0)$.
The set of rational points $X(F)$ decomposes into
a disjoint union

\[
C_1=X(F_0)=G\cdot P_1,\ \ \ \ \ \ 
C_2=X(F)-X(F_0)=G\cdot P_2,
\]
with $|C_1|=p+1$ and $|C_2|=2p(p-1)$.

\item[(b)] Case $p\equiv 1\, ({\rm mod}\, 4)$:

The automorphism group of $X/F$ acts transitively
on $X(F)$ and the stabilizer of any point is
a group of order $2p(p-1)$.

\end{itemize}
\end{proposition}

\begin{remark}
The proof of this proposition
is omitted, so may be regarded as a conjecture
instead, if the reader wishes.
It has been verified using MAGMA if $p=5, 7,11, 13$. 
It has been proven in an
email to the first author by Bob Guralnick.
\end{remark}

This and Lemma \ref{lemma:large} imply 
every point in $X(F)$ is ramified for the covering
$X\rightarrow X/G$.

Let $P_1=(1:0:1)$ and let $H_1$ be its stabilizer in $G$.
We have already seen that $H_1$ is a solvable group of order 
$2p(p-1)$ generated by 
$\gamma_1$, $\gamma_2(a)$, and $\gamma_3$. As a consequence,
$|C_1|=|G\cdot P_1|=|G|/|H_1|=p+1$ 

Using $H_1=\langle \gamma_1,\gamma_2,\gamma_3\rangle $ and the 
explicit expressions for the $\gamma_i$, it can be checked
directly that no $g\in H_1$, $g\not= 1$, fixes any
$P\in C_2$. Therefore, $H_1\cap H_2=\{1\}$.

According to the proposition, 
the stabilizer $H_2$ of $P_2$ has order $p+1$.
This and $|G|=|H_1|\cdot |H_2|$ implies $G=H_1\cdot H_2$.
In other words, $H_2$ is a complement of $H_1$ in $G$.
(As {\it sets}, $G/H_1\cong PGL_2(p)/B\cong \ppp^1(F_0)$.)

In fact,
if $B$ denotes the (Borel) upper-triangular subgroup
of $PGL(2,p)$ then $H_1=\tau^{-1}(B)$.
Since $B$ is solvable and any abelian cover of a solvable group is
solvable, $H_1$ is solvable.
Since $B$ is not normal in $PGL_2(p)$, $H_1$ is not normal in $G$.

By the proposition, the divisor $D_2$ associated to $O_2$ 
has degree $2p(-1)>2g$, so by the Riemann-Roch theorem, 
${\rm dim}(L(D_2))=2p(p-1)+1-\frac{p-1}{2}$. 
The theorem implies that, in this case, the largest 
irreducible constituent of $L(D_2)$ is dimension
$d_G=p$.

\section{Applications}
\label{sec:4}

In this section we discuss connections with the theory 
of error-correcting codes.

Throughout this section, we assume $X$, $G$, and $D$ are
as in Theorem \ref{thrm:main}. 
Assume $F$ is finite.

Let $P_1,...,P_n\in X(F)$ be distinct points
and $E=P_1+...+P_n\in Div(X)$ be stabilized by $G$. 
This implies that $G$ acts on the set
$supp(E)$ by permutation.
Assume ${\rm supp}(D)\cap {\rm supp}(E)=\emptyset$.
Let $C=C(D,E)$ denote the AG code

\begin{equation}
\label{eqn:AGcode}
C=\{(f(P_1),...,f(P_n))\ |\ f\in L(D)\}.
\end{equation}
This is the image of $L(D)$ under the evaluation map

\begin{equation}
\label{eqn:eval}
\begin{array}{c}
eval_E:L(D)\rightarrow F^n,\\
f \longmapsto (f(P_1),...,f(P_n)).
\end{array}
\end{equation}
The group $G$ acts on $C$ by $g\in G$ sending
$c=(f(P_1),...,f(P_n))\in C$ to $c'=(f(g^{-1}(P_1)),...,f(g^{-1}(P_n)))$,
where $f\in L(D)$.
First, we observe that this map, denoted $\phi(g)$, is well-defined.
In other words, if $eval_E$ is not injective and
$c$ is also represented by $f'\in L(D)$, so
$c=(f'(P_1),...,f'(P_n))\in C$, then we can easily verify
$(f(g^{-1}(P_1)),...,f(g^{-1}(P_n)))=
(f'(g^{-1}(P_1)),...,f'(g^{-1}(P_n)))$.
(Indeed, $G$ acts on the set
$supp(E)$ by permutation.)
This map $\phi(g)$ induces a homomorphism of $G$ into the permutation
automorphism group of the code ${\rm Aut}(C)$, denoted 

\begin{equation}
\label{eqn:phi}
\phi:G\rightarrow {\rm Aut}(C)
\end{equation}
(Prop. VII.3.3, \cite{Sti}, and \S 10.3, page 251, of 
\cite{St})\footnote{Both of these references define $\phi$ by 
$\phi(g)(c)=(f(g(P_1)),...,f(g(P_n)))$. However, this is a
homomorphism only when $G$ is abelian.}. 
The paper Wesemeyer \cite{W} investigated $\phi$ when
$C$ is a one-point AG code arising from a certain
family of planar curves.

\subsection{Separation of points}

To investigate the kernel of this map $\phi$, we introduce the following 
notion. Let $H\in Div(X)$ be any divisor.
We say that the space $L(H)$ 
{\bf separates points} if for
all points $P,Q\in X$,
$f(P)=f(Q)$ (for all $f\in L(H)$) 
implies $P=Q$ (see \cite{H}, chapter II, \S 7). 

We shall show that
Riemann-Roch spaces separate points for ``big enough'' divisors.

If $G$ is a group of automorphisms of $X$ defined over $F$
then $G$ induces an automorphism on the image of 
the evaluation map $eval_E : L(D)\rightarrow F^n$.
For this discussion, let us assume this is an injection.
(This is not a serious assumption.) 
To understand the kernel of this map $\phi$ in (\ref{eqn:phi}),
we'd like to know whether or not
$(f(P_1),...,f(P_n))=(f(g^{-1}P_1),...,f(g^{-1}P_n))$
implies $P_i=g^{-1}P_i$, for $1\leq i\leq n$.

Let $X$ be a plane curve with irreducible equation

\[
y^n + f_1(x)y^{n-1} + ... + f_{n-1}(x)y + f_n(x) = 0,
\]
where deg$(f_i(x))\leq i$, $1\leq i\leq n$. We assume $n\geq 2$
but we do not assume $X$ is non-singular.

Let $D$ be a divisor on $X$ and let $(x)_\infty$ be the point divisor of
$x$, so deg$(x)_\infty =n$.

Recall 
that the Riemann-Roch space $L(D)$ {\bf separates points}
if, for each pair $P,Q\in X -supp(D)$, $f(P)=f(Q)$
for all $f\in L(D)$ implies $P=Q$ \cite{H}.

\begin{lemma}
If $(x)_\infty \leq D$ then $L(D)$ separates points.
\end{lemma}

The hypothesis cannot be omitted.

\pf
Note that if $D'\leq D$ and $L(D')$ separates points
then $L(D)$ does too.

By hypothesis, $L((x)_\infty)\subset L(D)$.
By Proposition III.10.5 in \cite{Sti}, $x^iy^j\in L((x)_\infty)$,
for $0\leq j\leq n-1$ and $0\leq i\leq 1-j$. (Here
$0\leq i\leq 1-j$ means $i=0$ when $j\geq 1$.)
Let $P_1=(x_1,y_1)$ and $P_2=(x_2,y_2)$.
The condition $f(P_1)=f(P_2)$,
for all $f\in L((x)_\infty)$ implies 
$x_1^iy_1^j=x_2^iy_2^j$, for all $i,j$ as above.
In particular, we may take $(i,j)=(1,0)$
and $(i,j)=(0,1)$, so $P_1=P_2$. Therefore, 
$L((x)_\infty)$ separates points, and hence
$L(D)$ does too.
\qed

\vskip .1in

As the following example shows, the lemma is in some
sense best possible.

\vskip .1in

{\bf Example}:
Let $F=GF(9)$ and $X$ be the curve defined over $F$ by

\[
y^2=x^3-x.
\]
Let $P_\infty$ be the point at infinity on $X$. The
spaces $L(mP_\infty)$,
$1\leq m\leq 2$, do not separate points on $X$. Indeed, 
there are distinct points
$P,Q\in X(F)$ which have the same $x$-coordinate.
Since $L(2P_\infty)=\langle 1,x\rangle$, it cannot distinguish
them. On the other hand, by the Lemma, 
$L(3P_\infty)$ {\it must} separate points.
Indeed, $L(3P_\infty)=\langle 1,x,y\rangle$,
so from the reasoning in the above proof, it is obvious that
it does.

As a consequence of the lemma (changing variables if necessary), 
we see that, if for some
$a\in F$, $(x-a)_\infty \leq D$ then $L(D)$ separates points.

{\bf Question}: Is the converse also true?

\subsection{The kernel of $\phi$}

The paper by Wesemeyer \cite{W} investigated the
homomorphism $\phi:G\rightarrow {\rm Aut}(C)$ in
some special cases.
In general, if $L(D)$ separates points then
\[
{\rm Ker}(\phi)=\{g\in G\ |\ g(P_i)=P_i,\ 1\leq i\leq n\}.
\]
It is known (proof of Prop. VII3.3, \cite{Sti}) that
if $n>2g+2$ then $\{g\in G\ |\ g(P_i)=P_i,\ 1\leq i\leq n\}$
is trivial. Therefore, if $n>2g+2$ and 
$L(D)$ separates points then $\phi$ is injective.

\begin{example} 
Let $F=GF(7)$ and let $X$ denote the curve defined by

\[
y^2=x^7-x.
\]
This has genus $3$. The automorphism group
$Aut_F(X)$ is a central 2-fold cover of $PSL_2(F)$: 
we have a short exact sequence,

\[
1\rightarrow H \rightarrow Aut_F(X)\rightarrow PSL_2(7)\rightarrow 1,
\]
where $H$ denotes the subgroup of $Aut_F(X)$ 
generated by the hyperelliptic involution (which happens 
to also be the center of $Aut_F(X)$).
(Over the algebraic closure $\overline{F}$,
$Aut_{\overline{F}}(X)/center \cong PGL_2(\overline{F})$,
by \cite{G}, Theorem 1.)
Generators for the automorphism group are given in (\ref{eqn:3.2}) above,
taking $p=7$.

There are $8$ $F$-rational points\footnote{MAGMA views the 
curve as embedded in a weighted projective space, 
with weights $1$, $4$, and $1$, in which the point at 
infinity is nonsingular.}:

\[
X(F)=\{
P_1=(1:0:0),
P_2=(0:0:1),
P_3=(1:0:1), ...,
P_8=(6:0:1)\}.
\]
The automorphism group acts transitively on $X(F)$.
Consider the projection $C\rightarrow \ppp^1$ defined by
$\phi(x,y)=x$. The map $\phi$ is ramified at every point in
$X(F)$ and at no others.

Let $G=Stab(P_1,Aut_F(X))$ denote the stabilizer of
the point at infinity in $X(F)$. 
All the stabilizers $Stab(P_i,Aut_F(X))$
are conjugate to each other in $Aut_F(X)$, $1\leq i\leq 8$.
The group $G$ is a non-abelian group of
order $42$ (In fact, the group $G/Z(G)$ is the non-abelian group
of order $21$, where $Z(G)$ denotes the center of $G$.)

It is known (Proposition VI.4.1, \cite{Sti}) that, 
for each $m\geq 1$, the Riemann-Roch space
$L(mP_1)$ has a basis consisting of monomials,

\[
x^iy^j,\ \ \ \ \ \ 
0\leq i\leq 6,\ j\geq 0,\ 2i+7j\leq m.
\]
Let $D=5P_1$, $S=C(F)-\{P_1\}$, and let

\[
C(D,S)=\{(f(P_2), ...,f(P_8))\ |\ f\in L(D)\}.
\]
This is a $[7,3,5]$ code over $F$.
In fact, $\dim(L(D))=3$, so the evaluation map
$f\longmapsto (f(P_2), ...,f(P_8))$,
$f\in L(D)$, is injective.
Since $G$ fixes $D$ and preserves $S$, 
it acts on $C$ via 

\[
g:(f(P_2), ...,f(P_8))\longmapsto 
(f(g^{-1}P_2), ...,f(g^{-1}P_8)),
\]
for $g\in G$.

Let $P$ denote the permutation group of
this code. It a group of order $42$.
However, it is not isomorphic to $G$.
In fact, $P$ has trivial center.
The (permutation) action of $G$ on this code implies
that there is a homomorphism 

\[
\psi :G\rightarrow P.
\]
What is the kernel of this map?
There are two possibilities:
either a subgroup of order $6$ or a subgroup of order $21$
(this is obtained by matching possible orders of
quotients $G/N$ with possible orders of subgroups of $P$).
Take the automorphisms $\gamma_1$, $\gamma_2$ with $a=2$ and $\gamma_3$.
If we identify $S=\{P_2,...,P_8\}$ with 
$\{1,2,...,7\}$ then 

\[
\gamma_1 \leftrightarrow (2,7)(3,6)(4,5)=g_1,
\]
\[
\gamma_2 \leftrightarrow (2,5,3)(4,6,7)=g_2,
\]
\[
\gamma_3 \leftrightarrow (1,2,...7)=g_3.
\]
The group $ker(\phi)=N=\langle g_2,g_3\rangle$ 
is a non-abelian normal subgroup of
$G=\langle g_1,g_2,g_3\rangle$ of order $21$.

\end{example}

\subsection{Permutation representations}

In this subsection, we show how theorems about AG codes
can, in some cases, give theorems about representations on
Riemann-Roch spaces.

Assume that $X/\fff$ is a hyperelliptic curve defined 
over a finite field $\fff$ 
of characteristic $p>2$ with automorphism 
group $G=Aut_\fff(X)$. 
Let $D$ be a $G$-equivariant
divisor on $X$, let ${\cal O}\subset X(\fff)$ be a 
$G$-orbit disjoint from the support of $D$, and let
$E=\sum_{P\in {\cal O}}P$. Let $P$ be the permutation automorphism
group of the code $C=C(D,E)$ defined in 
(\ref{eqn:AGcode}). 

Theorem 4.6 in \cite{W} implies that
if $n={\rm deg}(E)$ and
$t={\rm deg}(D)$ satisfy $n>{\rm max}(2t,2g+2)$ then
the map $\phi:G\rightarrow P$ is an isomorphism.
Using this, we regard $C$ as a $G$-module.
In particular, the (bijective) evaluation map 
$eval_E:L(D)\rightarrow C$ in (\ref{eqn:eval})
is $G$-equivariant.
Since $G$ acts (via its isomorphism with $P$) as a permutation
on $C$, we have proven the following result.

\begin{proposition}
Under the conditions above, the representation $\rho$ of $G$ on $L(D)$
is equivalent to a representation $\rho'$ with with property that,
for all $g\in G$, $\rho'(g)$ is a permutation matrix.

\end{proposition}

\subsection{Memory application}

If $C$ is an linear code with non-trivial permutation group
then this extra symmetry of the code may be useful in practice.
In order to store the elements of $C$, we need only store
one element in each $G$-orbit, so this symmetry
can be used to more efficiently
store codewords in memory on a computer.

\begin{example}
Let $G=S_3$ act on the genus $3$ 
Fermat quartic $X$ whose projective model is 
$x^4+y^4+z^4=0$ over $\fff_9=\fff_3(i)$,
where $i$ is a root of the irreducible polynomial
$x^2+1\in \fff_3[x]$. One can check that there are
exactly $6$ distinct points in the $G$-orbit
of $[\alpha:1:0]\in X(\fff_9)$, where $\alpha$ is a generator
of $\fff_9^\times$. Let
\[
\begin{array}{c}
G\cdot  [\alpha :1:0]=\{Q_1,...,Q_6\},\\
\ \ \ E=Q_1+...+Q_6\in Div(X),\ \ D=6\cdot [1:1:1]\in Div(X).
\end{array}
\]
Then $L(D)$ is $4$-dimensional, by the Riemann-Roch theorem. 
Note that no $Q_i$ belongs to the support of $D$,
so we may construct the Goppa code
\[
C=\{(f(Q_1),...,f(Q_6))\ |\ f\in L(D)\},
\]
a generator matrix being given by the $4\times 6$ matrix
$M=(f_i(Q_j))_{1\leq i\leq 4,1\leq j\leq 6}$,
where $f_1,...,f_4$ are a basis of $L(D)$.
According to \cite{M}, $dim_{\fff_9}(C)=4$
and the minimum distance of $C$ is $2$.
The action of an element in the group $G$ on $C$ permutes the $Q_i$, hence 
may be realized by permuting the coordinates of
each codeword in $C$ in the obvious way.
(In other words, the action of $G$ on $C$ is isomorphic to
the regular representation of $S_3$ on itself.)
Using the group action, storing all $|C|=9^4$ elements may be
reduced to storing only the representatives of each orbit
$C/S_3$. 

\end{example}

\subsection{Permutation decoding application}

If $C$ is an linear code with non-trivial permutation group
then this extra symmetry of the code may be useful in 
decoding.
Permutation decoding is discussed, for example, in 
Huffman and Pless \cite{HP}. We recall briefly, for the
convenience of the reader, the main ideas.

We shall assume that $C$ is in standard form.
Let $C$ be a $[n,k,d]$ linear code over $GF(q)$,
let $t=[(d-1)/2]$, and let $G=(I_k,A)$ denote the generator matrix
in standard form. From this matrix $G$, it is well-known and easy to
show that one can compute
an encoder $E:GF(q)^k\rightarrow GF(q)^n$ with
image $C$, and a parity check matrix $H=(B,I_{n-k})$ 
in standard form, $B=-A^t$.

The key lemma is the following result: 
Suppose $v=c+e$, where $c\in C$ and $e\in GF(q)^n$
is an error vector with Hamming weight $wt(e)\leq t$. 
Under the above conditions,
the information symbols of $v$ are correct if and only if 
$wt(Hv)\leq t$. 

Let $P$ denote the permutation automorphism group of $C$. 
The permutation algorithm is:

\begin{enumerate}
\item
For each $p\in P$, compute $wt(H(pv))$ until one with
$wt(H(pv))\leq t$ is found (if none is found, the algorithm fails).

\item
Extract the information symbols from $pv$,
and use $E$ to compute codeword $c_p$ from them.

\item
Return $p^{-1}c_p={\rm Decode}(v)$. 
\end{enumerate}
For example, if $P$ acts transitively then 
permutation decoding will correct at least one error. 

The key problem is 
to find a set of permutations in $P$ which moves 
the non-zero positions in every possible error vector
of weight $\leq t$ out of the information positions.
(This set, called a {\it PD-set}, will be used in step
1 above instead of the entire set $P$.)

\begin{example}
We give two examples of MDS codes 
for which permutation decoding applies.

\begin{enumerate}
\item
This is an example of a $[7,3,5]$ one-point AG
code over $GF(7)$ arising from the hyperelliptic curve
$y^2=x^5-x$.

{\small{
\begin{verbatim}
p:=7;
F:=GF(p);
P<x>:=PolynomialRing(F);
f:=x^p-x;
X:=HyperellipticCurve(f);
Div := DivisorGroup(X);
Pls:=Places(X,1);
S:=[Pls[i] : i in [2..#Pls]];
m:=4;
D := m*(Div!Pls[1]);
AGC := AlgebraicGeometricCode(S, D);
Length(AGC);
Dimension(AGC);
MinimumDistance(AGC);
WeightDistribution(AGC);
PG := PermutationGroup(AGC);IdentifyGroup(PG);
ZP:=Center(PG);IdentifyGroup(PG/ZP);
IsTransitive(PG);
GeneratorMatrix(AGC);
\end{verbatim}
}}
This code has generator matrix in standard form given by

\[
G=
\left(
\begin{array}{ccccccc}
1 &0 &0 &2 &5 &1 &5\\
0 &1 &0 &1 &5 &5 &2\\
0 &0 &1 &5 &5 &2 &1
\end{array}
\right).
\]
Moreover, the permutation automorphism group of the code is a group
of order $42$ generated by

\[
S=\{ (1, 7)(2, 6)(3, 4), (1, 4, 5)(2, 6, 3),  (1, 3)(2, 4)(5, 6)\}.
\]
The elements of $S\, \cup \, S\cdot S$ can be used as a PD-set,
where $S\cdot S=\{s_1s_2\ |\ s_i\in S\}$.

\item
This is an example of a $[13,5,9]$ one-point AG
code over $GF(13)$ arising from the hyperelliptic curve
$y^2=x^{13}-x$. Similar MAGMA commands, but with $p=13$, 
yields that his code has generator matrix in standard form given by

\[
G=
\left(
\begin{array}{ccccccccccccc}
1 &  0 &  0 &  0 &  0 &  3 &  2 &  3 &  8 &  6 & 10 &  7 & 12\\
0 &  1 &  0 &  0 &  0 &  3 & 12 &  1 &  5 & 11 &  4 & 10  & 5\\
0 &  0 &  1 &  0 &  0 & 11 &  8  & 7 &  2  & 2 &  4 &  6 & 11\\
0 &  0 &  0 &  1 &  0 &  6 &  9 & 10 &  4 & 11 & 10 & 11  & 3\\
0 &  0 &  0 &  0 &  1 &  4 &  9 &  6 &  8 & 10 & 12 & 6 &  9
\end{array}
\right).
\]
Moreover, the permutation group is generated by

\[
\begin{array}{c}
\{p_1=(1, 2)(3, 8)(4, 12)(5, 7)(6, 9)(10, 11),\\
\ \ \ \ \ \ \ \ \ \ p_2=(2, 3, 4, 5, 6, 7, 8, 9, 10, 11, 12, 13)\}.
\end{array}
\]
We shall show that $P-\{1\}$ can be chosen to be
a PD-set for $\leq 4$ errors. The argument proceeds on a case-by-case basis.
One of the worst cases is when there are {\it errors 
in positions $1$, $2$, $12$ and $13$}. In this case, apply
(reading right-to-left) $p_1p_2p_1p_2^3$. This pushes the 
error from positions $(1,2,12,13)$ to $(8,11,13,6)$,
in particular out of the information positions\footnote{Additionally, 
this algorithm will, in {\it some}
cases, work even if there are $5$ errors (e.g., in positions $1,2,3,4,5$).
Because $d=9$, with $5$ errors we cannot be sure that the 
permutation decoded vector is the one which was sent.}.
\end{enumerate}

\end{example}

\begin{conjecture} For one-point AG codes $C$ associated
to $y^2=x^p-x$ over $GF(p)$ of length $n=p$,
permutation decoding always applies.
Its complexity is
at worst the size of the permutation group of $C$, 
which we conjecture to be $O(p^2)=O(g^2)=O(n^2)$.
\end{conjecture}

This matches the complexity of some known algorithms.

\begin{conjecture}
For one-point AG codes $C$ associated
to $y^2=x^p-x$ over $GF(p^2)$ of length $n=2p(p-1)+p$,
permutation decoding always applies and is more efficient
in terms of computational complexity that the standard
decoding algorithm in \cite{St}. We conjecture that, if the points
in $X(F)$ are arranged suitably then the image of the
$Aut_F(X)$ in the permutation group of $C$ may be used
as a PD-set. Its complexity is
at worst the size of the automorphism group of $X$, 
which is $O(p^2)=O(g^2)=O(n)$.
\end{conjecture}

If true, to our knowledge, this beats the complexity of other
decoding algorithms, such as those in \cite{Sti}.

\begin{example}
We give examples of two AG codes 
for which permutation decoding probably applies.

\begin{itemize}
\item
This is an example of a $[91,5,66]$ code constructed from 
the trace of a $[91,3,87]$ one-point AG
code over $F=GF(49)$ arising from the hyperelliptic curve
$y^2=x^7-x$. (We use the trace code only because MAGMA
version 2.10 cannot compute the permutation group of a
code over $GF(49)$.)  

{\small{
\begin{verbatim}
p:=7;
F:=GF(p^2);
P<x>:=PolynomialRing(F);
f:=x^p-x;
X:=HyperellipticCurve(f);
Div := DivisorGroup(X);
Pls:=Places(X,1);
S:=[Pls[i] : i in [2..#Pls]];
m:=4;
D := m*(Div!Pls[1]);
AGC := AlgebraicGeometricCode(S, D);
Length(AGC);
Dimension(AGC);
MinimumDistance(AGC);
WeightDistribution(AGC);
AGC0:=Trace(AGC,GF(p));
Length(AGC0);
Dimension(AGC0);
MinimumDistance(AGC0);
WeightDistribution(AGC0);
PG := PermutationGroup(AGC0);
#PG; 
ZP:=Center(PG);
#ZP; 
IsTransitive(PG);
GeneratorMatrix(AGC0);
\end{verbatim}
}}
The permutation group of the trace of the AG code is 
huge: $1073852196$ elements. The 
automorphism group $G=Aut_F(X)$, which has $672$ elements, 
of the curve acts on $X(F)$ with only two orbits,
$O_1=X(GF(7))$ of size $8$ and $O_2=X(F)-O_1$ of
size $84$. (This follows from Proposition \ref{prop:orbits}
but was verified using MAGMA in this case.)
Of course, in a practical application,
one would want to index the points of
$X(F)$ so that the information positions are contained in 
$O_2$.

\item
Let $E$ denote the sum of all the points in $O_2$
and let $D$ be the sum of all the points in $O_1$.
Note $E$ has degree $84$, $D$ degree $8$, and $X$ has genus $3$.
Therefore, by Theorem 4.6 in \cite{W},
the permutation automorphism group $P$ of the AG code $C=C(D,E)$  
satisfies $P\cong G$. In other words, the map
$\phi$ in (\ref{eqn:phi}) is injective (which also
follows from the discussion above) and surjective.
\end{itemize}
\end{example}

\vskip .2in
{\it Acknowledgements}: We thank Keith Pardue,
Niels Borne, Bob Guralnick, Derek Holt 
and Amy Ksir, for kindly answering our 
questions on group theory and algebraic geometry. We thank 
Nick Sheppard-Barron for the reference to \cite{MP}.
Finally, we especially thank Bernhard K\"ock for 
many detailed comments improving the content of the
first version of this paper and for the references \cite{Ka}, \cite{K}.

\end{document}